\input amstex
\documentstyle{amsppt}
\topmatter
\title{Elementary Proof of MacMahon's Conjecture}\endtitle
\author{David M.\ Bressoud}\endauthor
\date{February 6, 1997}\enddate
\thanks Copyright to this work is retained by the
author. Permission is granted for the noncommercial
       reproduction of the complete work for educational or research purposes, and for the use of figures,
       tables and short quotes from this work in other books or journals, provided a full bibliographic
       citation is given to the original source of the material. \endthanks
\abstract  Major Percy A.\ MacMahon's first paper on plane partitions
\cite{4} included a conjectured generating function for symmetric plane
partitions. This conjecture was proven almost simultaneously by George Andrews and
Ian Macdonald, Andrews using the machinery of basic hypergeometric series
\cite{1} and Macdonald employing his knowledge of symmetric functions
\cite{3}. The purpose of this paper is to simplify Macdonald's proof by
providing a direct, inductive proof of his formula which expresses the sum of Schur
functions whose partitions fit inside a rectangular box as a ratio of determinants.
\endabstract
\address  Dept.\ of Mathematics \& Computer Science, Macalester College,
St.\ Paul, MN 55105, USA 
\endaddress
\email bressoud\@macalester.edu \endemail
\endtopmatter
\define\pf{\noindent {\bf Proof:} \ }
\define\inv{{\Cal{I}}}
\define\rhs{\text{RHS}}
\document

\noindent By a {\bf plane partition}, we mean a finite set, $\Cal{P}$, of lattice
points with positive integer coefficients, $\{(i,j,k)\} \subseteq {\Bbb N}^3$, with
the property that if $(r,s,t) \in \Cal{P}$ and $1 \leq i \leq r,\ 1 \leq j \leq s,\ 1
\leq k
\leq t$, then $(i,j,k)$ must also be in $\Cal{P}$. A plane partition is {\bf
symmetric} if $(i,j,k) \in \Cal{P}$ if and only if $(j,i,k)\in \Cal{P}$. MacMahon's
conjecture states that the generating function for symmetric plane partitions whose
$x$ and $y$ coordinates are less than or equal to $n$ and whose $z$ coordinate is
less than or equal to $m$ is given by
$$
\prod_{i=1}^n \frac{1-q^{m+2i-1}}{1-q^{2i-1}} \prod_{1 \leq i < j \leq
n}\frac{1-q^{2(m+i+j-1)}}{1-q^{2(i+j-1)}}.
$$

Our proof parallels that of Ian Macdonald \cite{3} which divides
into three distinct pieces. We shall concentrate on the middle piece which is
the most difficult and the heart of his argument. Macdonald derived it as a
corollary of a formula for Hall-Littlewood polynomials. Details of the proof
of Macdonald's formula as well as a generalization may be found in \cite{2}. We shall
prove the middle piece directly by induction on the number of variables.

The first piece of Macdonald's proof is the observation, known before Macdonald,
that there is a one-to-one correspondence, preserving the number of lattice points,
between bounded symmetric plane partitions and {\bf column-strict plane partitions}
with $y$ coordinates bounded by $m$, $z$ coordinates bounded by $2n-1$, and in which
and non-empty columns have odd height. The column at position $(i,j)$ is the set
of $(i,j,k) \in \Cal{P}$, and the column height is the cardinality of this set. To
say that the partition is column-strict means that if $1\leq h < i$ and the column
at $(h,j)$ is non-empty, then the column height at $(h,j)$ must be strictly greater
than the column height at $(i,j)$.

From this observation and the definition of the Schur function, $s_{\lambda}$, as a
sum over semi-standard tableaux of shape $\lambda$, it follows that the generating
function for bounded symmetric plane partitions is given by
$$ \sum_{\lambda\subseteq \{m^n\}}
s_{\lambda}(q^{2n-1},q^{2n-3},\ldots,q), $$
where the sum is over all partitions, $\lambda$, into at most $n$ parts each of
which is less than or equal to $m$.

The second piece of Macdonald's proof is the following theorem which is the result
that we shall prove in this paper. 

\medskip

\noindent{\bf Theorem} {\it For arbitrary positive integers $m$ and $n$,\/}
$$
\sum_{\lambda\subseteq \{m^n\}} s_{\lambda}(x_1,\ldots,x_n) =
\frac{\det(x_i^{j-1}-x_i^{m+2n-j})}{\det(x_i^{j-1}-x_i^{2n-j})}.\tag1
$$

The final piece of Macdonald's proof is to rewrite the right side of 
equation~(1) when $x_i = q^{2(n-i)+1},\ 1\leq i \leq n,$ as a ratio of
products by employing the Weyl denominator formula for the root system $B_n$:
$$ \det(x_i^{j-1}-x_i^{2n-j}) = \prod_{i=1}^n (1-x_i) \prod_{1\leq i <
j \leq n} (x_i-x_j)(x_i\, x_j-1). \tag2 $$

There is a very simple inductive proof of this case of the Weyl denominator
formula. Let  $D_n(x_1,...x_n) = \det( x_j^{i-1} - x_j^{2n-i} )$. This is a
polynomial of degree $2n-1$ in $x_1$ with roots at $1, x_2, .., x_n, x_2^{-1}, ...,
x_n^{-1}$. The coefficient of $x_1^{2n-1}$ is $- x_2 \cdots x_n
D_{n-1}(x_2,...,x_n)$.

Before we begin the proof of the theorem, we note that it similarly implies
Gordon's identity (\cite{3}, page 86):
$$ \sum_{\lambda\subseteq \{m^n\}} s_{\lambda}(q^n,q^{n-1},\ldots,q) = \prod_{1\leq
i \leq j \leq n} \frac{1-q^{m+i+j-1}}{1-q^{i+j-1}}. $$

\head Proof of the theorem \endhead

\noindent We shall need the following lemma.

\medskip

\noindent{\bf Lemma}
$$\multline x_1 \cdots x_n \sum_{k=1}^n (-1)^{k-1}(1-x_k) x_k^{-1} \prod_{i
\neq k} (1-x_i\, x_k) \prod_{1\leq i < j \leq n \atop i,j \neq k}(x_j - x_i)
 \\
= (1 - x_1 \cdots x_n)\prod_{1\leq i < j \leq n}(x_j-x_i).   
\endmultline\tag3 $$

\pf We verify that this lemma is correct for $n=2$ or 3 and proceed by induction.
The left side of equation~(3) is an anti-symmetric polynomial. If we
divide it by
$\prod_{1\leq i < j
\leq n}(x_j-x_i)$, we obtain a symmetric polynomial. Let us
denote this ratio by 
$$ F(x_1,\ldots,x_n) = x_1 \cdots x_n \sum_{k=1}^n (1-x_k) x_k^{-1} 
\prod_{i \neq k} {1-x_i\, x_k \over x_i - x_k}. $$
As a function of $x_1$, $F$ is a polynomial of degree at most $n$ divided by a
polynomial of degree $n-1$, and is therefore a linear polynomial in $x_1$. It is
easily verified that
$$\align
F(0,x_2,\ldots,x_n) & =  1, \\
F(1,x_2,\ldots,x_n) & =  F(x_2,\ldots,x_n) \\
& =  1-x_2 x_3\cdots x_n. \qed
\endalign $$
   
We use equation~(2) to rewrite the right hand side of the
theorem as
$$ \frac{\det(x_i^{j-1}-x_i^{m+2n-j})}{\prod_{i=1}^n (1-x_i) \prod_{1\leq i <
j \leq n} (x_i-x_j)(x_i\, x_j-1)}. $$
We shall also use the representation of the Schur function as a ratio of
determinants:
$$ s_{\lambda}(x_1,\ldots,x_n) = \frac{\det(x_i^{\lambda_j+n-i})}{\prod_{1\leq i <
j \leq n}(x_i-x_j)}. $$
Combining these, our theorem can be restated as
$$
\det(x_i^{j-1}-x_i^{m+2n-j}) =
\sum_{\lambda\subseteq\{m^n\}}\det(x_i^{\lambda_j+n-j})\
\prod_{i=1}^n(1-x_i)\,\prod_{1\leq i < j \leq n} (x_i\, x_j -1). \tag4
$$
When we expand these determinants, we see that the theorem to be proved is
equivalent to
$$ \multline \sum_{\sigma,S} (-1)^{\inv(\sigma)+|S|} \prod_{i\in S}
x_i^{m+2n-\sigma(i)}\,\prod_{i\not\in S} x_i^{\sigma(i)-1}  \\
= \sum_{\lambda,\sigma} (-1)^{\inv(\sigma)} \prod_{i=1}^n x_i^{\lambda_{\sigma(i)}
+ n - \sigma(i)}\
\prod_{i=1}^n(1-x_i)\,\prod_{1\leq i < j \leq n} (x_i\, x_j -1),
\endmultline\tag5
$$
where $\inv(\sigma)$ is the inversion number. The first sum is over all
permutations, $\sigma$, and subsets, $S$, of $\{1,\ldots,n\}$. The second sum is
over partitions $\lambda \subseteq \{m^n\}$ and permutations.

Our proof will be by induction on $n$. It is easy to check that this equation is
correct for $n=1$ or 2. Let RHS denote the right hand side of
equation~(5). We shall sum over all possible values of $\lambda_n$ and
$k = \sigma^{-1}(n)$. Given $\lambda_n$ and $k$, we subtract $\lambda_n$ from
each part in $\lambda$ to get $\lambda' \subseteq \{(m-\lambda_n)^{n-1}\}$. The
permutation $\sigma$ is uniquely determined by $k$ and a one-to-one mapping
$\sigma':\{1,\ldots,n\}\backslash\{k\} \to \{1,\ldots,n-1\}$. We can express the
right hand side of equation~(5) as:
$$\split  \rhs & =  \sum_{\lambda_n=0}^m \sum_{k=1}^n
(-1)^{n+k}(1-x_k)x_k^{-1}(x_1\cdots x_n)^{\lambda_n+1} \prod_{i\neq k}(x_i\,x_k-1) \\
& \qquad \sum_{\lambda',\sigma'} (-1)^{\inv(\sigma')} \prod_{i\neq
k}x_i^{\lambda'_{\sigma'(i)} + (n-1) - \sigma'(i)}
\prod_{i=1 \atop i \neq k}^n(1-x_i)\,\prod_{1\leq i < j \leq n \atop i,j \neq k} (x_i
\,x_j -1).
\endsplit $$
We apply the induction hypothesis to the inner sum and then sum over $\lambda_n$:
$$\split
\rhs & =  \sum_{k=1}^n (-1)^{n+k}(1-x_k)\prod_{i\neq k}(x_i\, x_k-1) \\
& \qquad \sum_{\sigma,S} (-1)^{\inv(\sigma)+|S|} \prod_{i\in
S}x_i^{m+1+2n-2-\sigma(i)}\prod_{i\in\overline{S}}x_i^{\sigma(i) }
\frac{1-x_k^{m+1}\prod_{i\in\overline{S}}x_i^{m+1}}{1
- x_k\prod_{i\in\overline{S}}x_i},
\endsplit $$
where the inner sum is over all one-to-one mappings $\sigma$ from
$\{1,\ldots,n\}\backslash\{k\} \to \{1,\ldots,n-1\}$ and subsets $S$ of
$\{1,\ldots,n\}\backslash\{k\}$. We use $\overline{S}$ to denote the complement of
$S$ in $\{1,\ldots,n\}\backslash\{k\}$.

It is convenient at this point to replace $x_i^{m+1}$ by $t_i\, x_i^{2-2n}$ on each
side of the equation to be proved. Our theorem is now seen to be equivalent to
$$ \multline \sum_{\sigma,S} (-1)^{\inv(\sigma)+|S|} \prod_{i\in S}
t_i\, x_i^{1-\sigma(i)}\,\prod_{i\not\in S} x_i^{\sigma(i)-1}  \\
= {\split & \sum_{k=1}^n (-1)^{n+k}(1-x_k)\prod_{i\neq k}(x_i\, x_k-1)  \\
   & \qquad \sum_{\sigma,S} (-1)^{\inv(\sigma)+|S|} \prod_{i\in
S}t_i\, x_i^{-\sigma(i)}\prod_{i\in\overline{S}}x_i^{\sigma(i) }
\ \frac{1-\prod_{i\not\in S}t_i\,x_i^{2-2n}}{1
- \prod_{i\not\in S}x_i}.\endsplit} 
\endmultline\tag6 $$

The sum on $\sigma$ on the right hand side is a Vandermonde determinant in $n-1$
variables. We replace it with the appropriate product and then interchange the
summation on $S$, which must be a proper subset of
$\{1,\ldots,n\}$, and $k$, which cannot be an element of $S$: 
$$ \multline
\rhs  =  \sum_{S \subset \{1,\ldots,n\}}(-1)^{|S|}
\prod_{i\in S}t_ix_i^{-1} \prod_{i\not\in S}x_i \left(\frac{1-\prod_{i\not\in
S}t_i\,x_i^{2-2n}}{1 - \prod_{i\not\in S}x_i}\right) \\
 \sum_{k \not\in S} (-1)^{n+k}(1-x_k)x_k^{-1}\prod_{i\neq k}(x_i\,x_k
- 1)
\prod_{i<j \atop i,j\neq k}(x_j^{\epsilon_j}-x_i^{\epsilon_i}),
\endmultline $$
where $\epsilon_i = -1$ if $i\in S$, $= +1$ if $i \not\in S$. We rewrite
$$ \align (-1)^{n+k} \prod_{i\neq k}(x_i\,x_k-1)
& =  \prod_{i<k}(x_i\,x_k-1)\prod_{i>k}(1-x_i\,x_k) \\
& =  \prod_{i<k \atop i \not\in S}(x_i\,x_k-1)\prod_{i>k \atop i \not\in
S}(1-x_i\,x_k) \\
& \qquad \times\ \prod_{i\in S}x_i \prod_{i<k
\atop i \in S}(x_k-x_i^{-1})\prod_{i>k\atop i \in S}(x_i^{-1}-x_k), 
\endalign $$
and then factor all terms that involve $x_i$, $i\in S$, out of the sum on $k$. The
sum on $k\not\in S$ can now be evaluated using the lemma:
$$ \align
\rhs & =   \sum_{S \subset \{1,\ldots,n\}}(-1)^{|S|}\prod_{i\in S}
t_i \left(1-\prod_{i\not\in S}t_i\,x_i^{2-2n}\right)  \prod_{1\leq i < j \leq
n} (x_j^{\epsilon_j}-x_i^{\epsilon_i}) \\ & = 
\sum_{S,\sigma}(-1)^{\inv(\sigma)+|S|}\prod_{i\in S}t_ix_i^{1-\sigma(i)}
\prod_{i\not\in S}x_i^{\sigma(i)-1} \\
&   \qquad -\  t_1\cdots t_n \sum_{S,\sigma}(-1)^{\inv(\sigma)+|S|}\prod_{i\in
S}x_i^{1-\sigma(i)} \prod_{i\not\in S}x_i^{\sigma(i)+1-2n},
\endalign $$
where both sums are over all {\it proper\/} subsets $S$ of $\{1,\ldots,n\}$.
Equation~(6)---which we have seen is equivalent to the theorem---now
follows from the observation that when we sum over all subsets $S$ of
$\{1,\ldots,n\}$,
	$$ \sum_{S,\sigma}(-1)^{\inv(\sigma)+|S|}\prod_{i\in
S}x_i^{1-\sigma(i)} \prod_{i\not\in S}x_i^{\sigma(i)+1-2n} = \det(x_i^{j+1-2n} -
x_i^{1-j}) = 0. $$

\enddocument